\documentclass[a4paper]{article}

\usepackage{graphicx}
\usepackage{mathrsfs}
\usepackage[latin2]{inputenc}
\usepackage[magyar,english]{babel}
\usepackage{amsmath,amssymb,amscd}
\newtheorem{definicio}{Definíció}

\newtheorem{theorem}[definicio]{Theorem}

\makeatletter
\def\egyfill{$\m@th\hbox{\raisebox{-1.5pt}[1.5pt][0pt]{$\mathord\leftarrow\!$}}
\mkern-5mu
\cleaders\hbox{\raisebox{-1.5pt}[1.5pt][0pt]{$\!\mathord-\!$}}\hfill
\mkern-5mu
\mathord{\hbox{\raisebox{-1.5pt}[1.5pt][0pt]{$\!\mathord\rightarrow$}}}$}
\makeatother
\def\egy { \overset{\hbox{\egyfill}} }

\makeatletter
\def\fegyfill{$\m@th\hbox{\raisebox{-1.5pt}[1.5pt][0pt]{$\mathord-\!$}}
\mkern-5mu
\cleaders\hbox{\raisebox{-1.5pt}[1.5pt][0pt]{$\!\mathord-\!$}}\hfill
\mkern-5mu
\mathord{\hbox{\raisebox{-1.5pt}[1.5pt][0pt]{$\!\mathord\rightarrow$}}}$}
\makeatother

\author{Zoltán Szilasi}
\title{Hagge configurations and a projective generalization of inversion}
\date{}

\begin{document}

\maketitle

\section{Introduction}

The discovery of Hagge's circle by K. Hagge in 1907 \cite{Hagge} opened new perspectives in classical geometry (\cite{BradSmith}, \cite{J}, \cite{Preiser}, \dots). In a recent paper \cite{B}, Bradley described two generalizations of Hagge's theorems. By means of coordinate calculations, first he proved

\begin{theorem}
Let a triangle $ABC$ be given in the Euclidean plane. Let $D$ be a point, not lying on the side lines of the triangle, and let $\Sigma$ be a circle passing through $D$. If $\Sigma\backslash\left\{D\right\}$ meets the circles $BCD$, $ACD$, $ABD$ at the points $U$, $V$, $W$, and meets the lines $\egy{AD}$, $\egy{BD}$, $\egy{CD}$ at the points $X$, $Y$, $Z$, respectively, then the lines $\egy{UX}$, $\egy{VY}$, $\egy{WZ}$ are concurrent.
\end{theorem}

\begin{figure}[ht]
	\centering
			\includegraphics[scale=0.25]{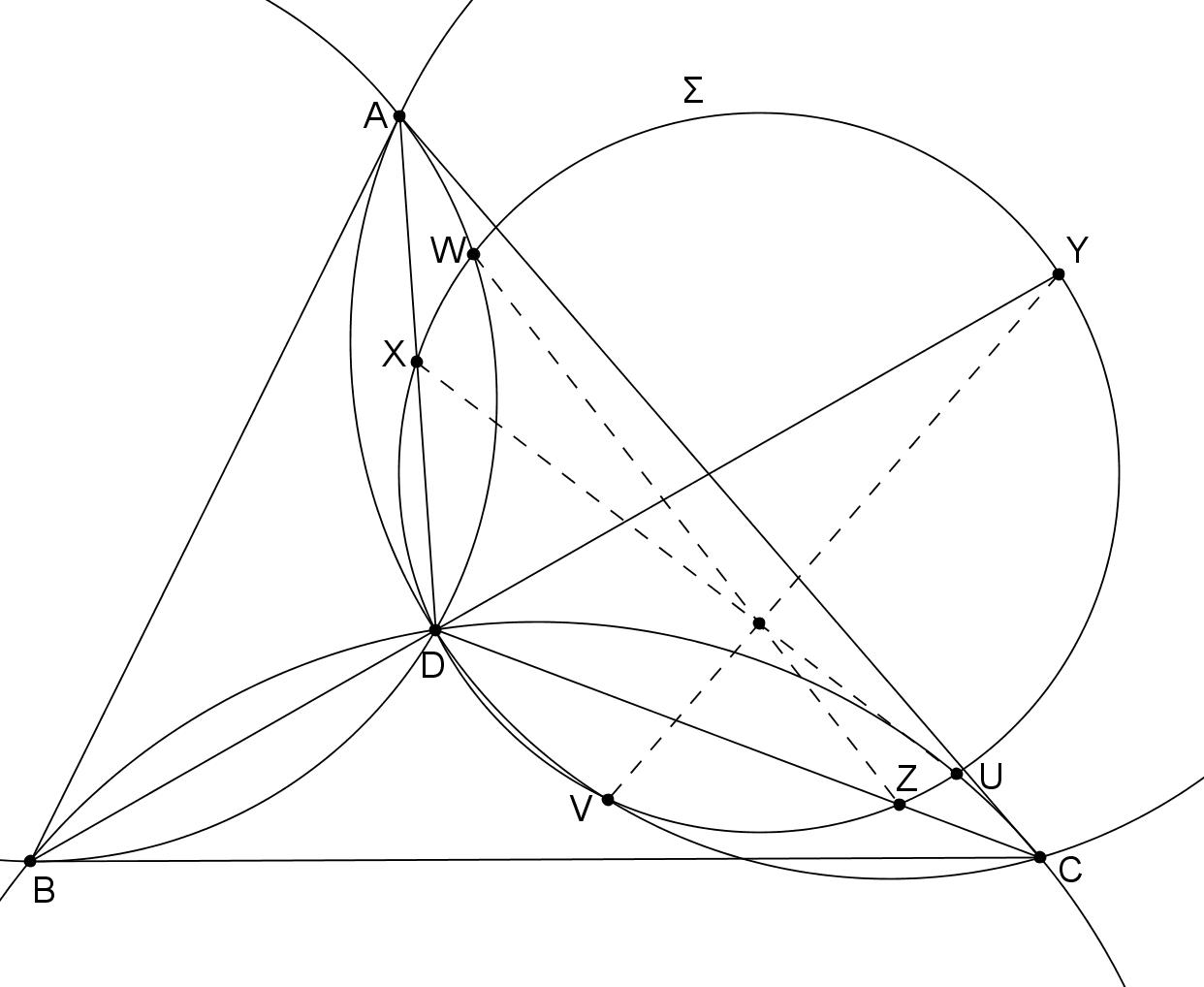}
\end{figure}

Next he deduced an essentially projective generalization of this result:

\begin{theorem}
Let a triangle $ABC$ be given in the Euclidean plane. Let $D$, $E$, $F$ be non-collinear points, neither of which lies on a side line of the triangle. If a conic $\Sigma$ passes through the points $D$, $E$, $F$, and $\Sigma\backslash\left\{D\right\}$ meets the conics $BCDEF$, $ACDEF$, $ABDEF$ at the points $U$, $V$, $W$, and meets the lines $\egy{AD}$, $\egy{BD}$, $\egy{CD}$ at the points $X$, $Y$, $Z$, then the lines $\egy{UX}$, $\egy{VY}$, $\egy{WZ}$ are concurrent.
\end{theorem}

\begin{figure}[ht]
	\centering
			\includegraphics[scale=0.25]{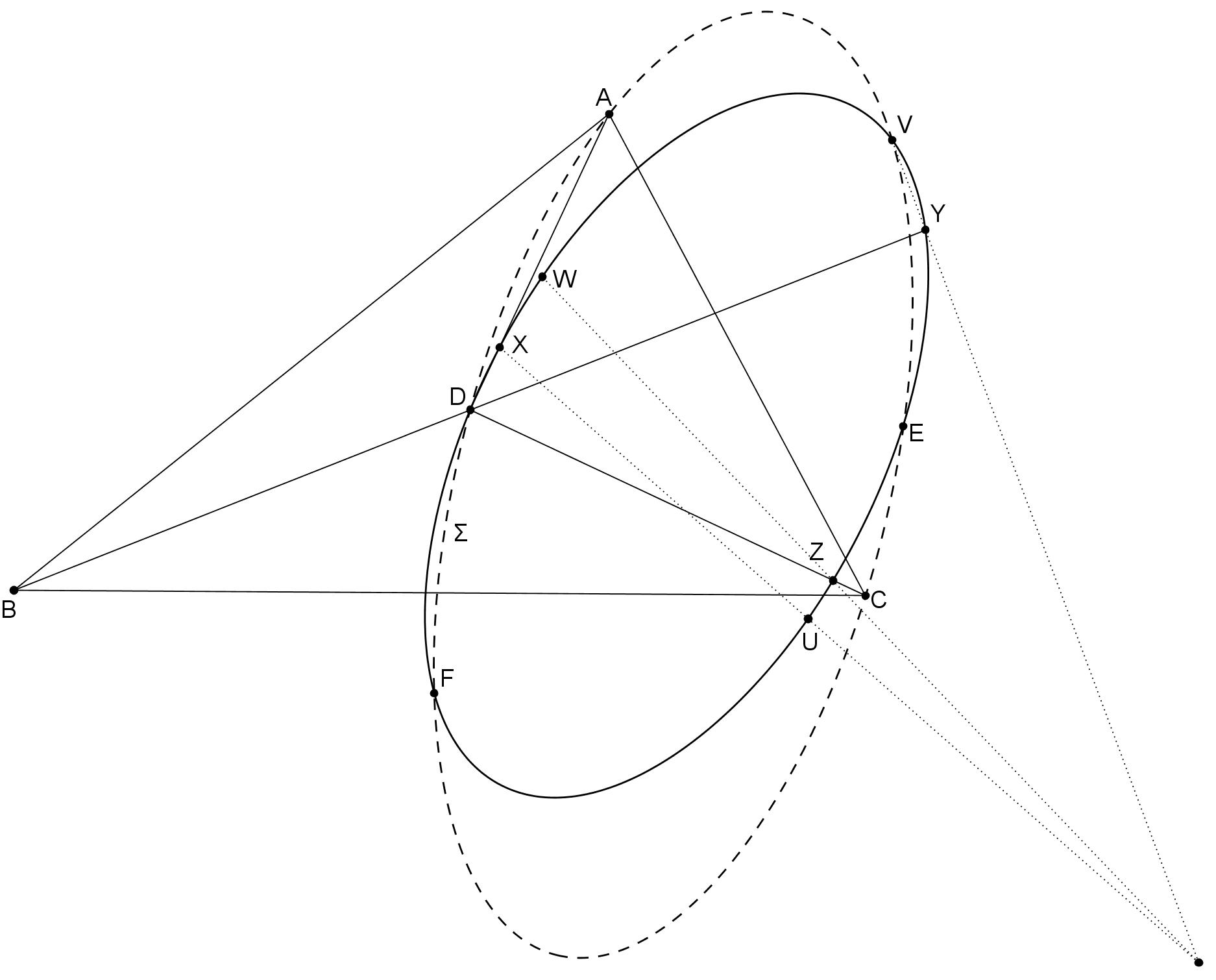}
\end{figure}

Theorem 2 indeed reduces to Theorem 1 if $E$ and $F$ are the 'circular points at infinity'. In this note we present synthetic, elementary proofs for both of these theorems. Our proof for Theorem 2 does not rely on Theorem 1 (so we may immediately deduce the first theorem from the second one). The reasoning applied in the proof of Theorem 2 is a substantial refinement of that in the proof of Theorem 1. In fact we show that Theorem 2 is valid in any Pappian projective plane satisfying Fano's axiom.

In both proofs we need the following basic facts from projective geometry.\\

\emph{Let a Pappian plane be given, satisfying Fano's axiom. Then we have}\\

\textbf{(A)} \emph{The three pairs of opposite sides of a complete quadrangle meet any line (not passing through a vertex) in the three pairs of an involution.}\\

\textbf{(B)} \emph{If U, V, W, X, Y, Z are six points on a conic, then the three lines} $\egy{UX}$, $\egy{VY}$, $\egy{WZ}$ \emph{are concurrent, if and only if, (UX), (VY), (WZ) are pairs of an involution on the conic.}\\

For a proof we refer to Coxeter's book \cite{Cox1}.

\section{An elementary proof of Theorem 1}

We may interpret the Euclidean plane as a part of its projective closure. The latter is a Pappian plane satisfying Fano's axiom (in fact, it is isomorphic to the real projective plane).

Apply an inversion of pole $D$, denoting the images of points and sets by a prime. Then the sets $$\left\{D,A',X'\right\}\textrm{  ,  }\left\{D,B',Y'\right\}\textrm{  ,  }\left\{D,C',Z'\right\}\textrm{  ,  }\left\{A',B',W'\right\}\textrm{  ,  }\left\{A',C',V'\right\}\textrm{  ,  }\left\{B',C',U'\right\}$$ are collinear, therefore the opposite sides of the complete quadrangle $A'B'C'D$ meet the line $\Sigma'$ at the pairs $(U'X')$, $(V'Y')$, $(W'Z')$. By \textbf{(A)}, these are the three pairs of an involution. On the other hand, inversion preserves cross ratio cr, and pairs  $(P_1P'_1)$, $(P_2P'_2)$, $(P_3P'_3)$ are pairs of an involution, if and only if, $\textrm{cr}(P_1,P_2,P_3,P'_3)=\textrm{cr}(P'_1,P'_2,P'_3,P_3)$. Thus the involution sends the pairs $(U'X')$, $(V'Y')$, $(W'Z')$ to the pairs of an involution on the circle $\Sigma=(\Sigma')'$. Hence, by \textbf{(B)}, the lines $$\egy{UX}\textrm{  ,  }\egy{VY}\textrm{  ,  }\egy{WZ}$$ are concurrent.

\begin{figure}[ht]
	\centering
			\includegraphics[scale=0.20]{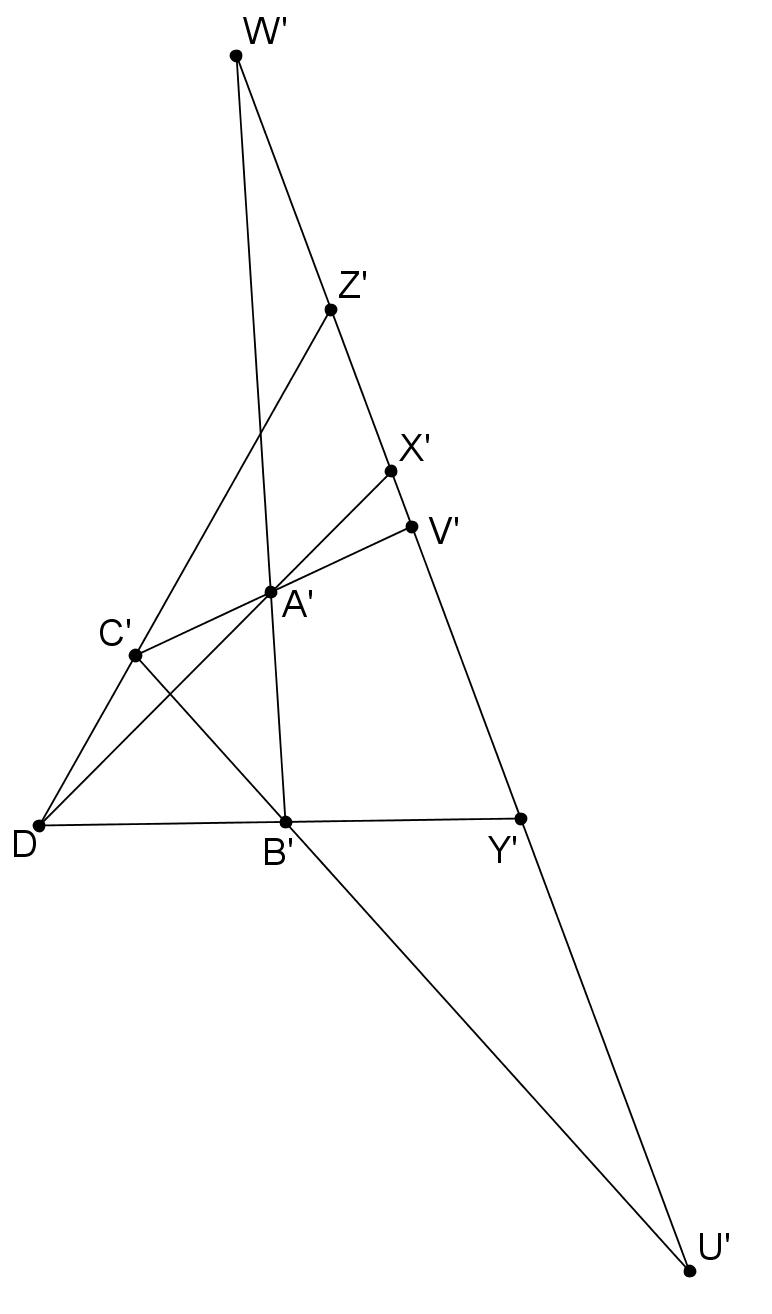}
\end{figure}

\section{An elementary proof of Theorem 2}

In this section we consider a Pappian projective plane $\mathbf{P}$ satisfying Fano's axiom. First we collect some basic facts concerning the so-called \emph{Steiner-correspondence}. These suggest that it is a good candidate for a purely projective generalization of inversion. Indeed, Steiner correspondence will play the same role in the proof of Theorem 2 as inversion in the proof of Theorem 1.

Let $\Sigma_1$ and $\Sigma_2$ be two fixed conics in $\mathbf{P}$.\\

(1) It is known (see e.g. \cite{Hatt}) that $\Sigma_1$ and $\Sigma_2$ have a common self-polar triangle $\Delta=\left\{D_1,D_2,D_3\right\}$. For every point $P\in\mathbf{P}\backslash\Delta$ let $p_1$ and $p_2$ be the polars of $P$ with respect to $\Sigma_1$ and $\Sigma_2$, respectively. Then the mapping $$\mathcal{S}: \mathbf{P}\backslash\Delta\rightarrow\mathbf{P}\textrm{  ,  }P\mapsto P':=p_1\cap p_2$$ is said to be the \emph{Steiner correspondence} with respect to $\Sigma_1$ and $\Sigma_2$. If $P'=\mathcal{S}(P)$, then we say that the points $P$ and $P'$ are in \emph{Steiner-correspondence}. Notice that if $P$ is a vertex of $\Delta$, e.g., $P=D_1$, then $p_1=p_2=\egy{D_2D_3}$, so the Steiner correspondence cannot be defined. Over $$\mathbf{P}\backslash\left\{\egy{D_1D_2}\cup\egy{D_2D_3}\cup\egy{D_3D_1}\right\}$$ $\mathcal{S}$ is involutive, and hence invertible.\\

(2) Let $l\subset\mathbf{P}$ be a line, \emph{not} passing through any vertices of $\Delta$. We show that the set $$\mathcal{S}(l)=\left\{\mathcal{S}(P)\in\mathbf{P}\textrm{  }\vert\textrm{  }P\in l\right\}$$ is a \emph{conic}. Let $L_1$ and $L_2$ be the poles of $l$ with respect to $\Sigma_1$ and $\Sigma_2$, and consider the pencils $\mathcal{L}_i$ with centres $L_i$, $i\in\left\{1,2\right\}$. Then for each point $P\in l$ the point $P'=\mathcal{S}(P)$ can be obtained as the intersection of two corresponding lines in $\mathcal{L}_1$ and $\mathcal{L}_2$. Since there is a projectivity between $\mathcal{L}_i$ and the range of all points on $l$ for $i\in\left\{1,2\right\}$, it follows that we also have a projectivity $f:\mathcal{L}_1\rightarrow\mathcal{L}_2$. Then, by Steiner's characterization of conics, the locus of points $m\cap f(m)\textrm{ , }m\in\mathcal{L}_1$ is a conic. Clearly, this conic is just the set $\mathcal{S}(l)$.\\

(3) Observe that the conic $\mathcal{S}(l)$ contains the vertices of $\Delta$, since $\mathcal{S}$ sends every side line of $\Delta$ into the vertex opposite to the side. From the involutiveness of $\mathcal{S}$ it follows that \emph{the image of a conic passing through a vertex of }$\Delta$\emph{ is a line}.\\

(4) By the reasoning applied in (2) we can also see that the Steiner correspondence sends the pairs of an involution of points on $l$ to the pairs of an involuion of points on the conic $\mathcal{S}(l)$.\\

(5) Suppose, finally, that the line $l\subset\mathbf{P}$ \emph{passes through} a vertex $D\in\Delta$. We claim that in this case the image of $l\backslash\left\{D\right\}$ under $\mathcal{S}$ is a range of points on a line. Indeed, using the same notation as in (2), the poles $L_1$ and $L_2$ are on the side line $d$ opposite to $D$. $d$ is the polar of $D$, so the projectivity $f:\mathcal{L}_1\rightarrow\mathcal{L}_2$ sends $d$ into itself, therefore $f$ is a perspectivity, and the points $m\cap f(m)$ $(m\in\mathcal{L}_1)$ are collinear. Again, this point set is just $\mathcal{S}(l\backslash\left\{D\right\})$.\\

Now we are in a position to prove Theorem 2 in the given Pappian plane $\mathbf{P}$. Consider two conics $\Sigma_1$, $\Sigma_2$ with the same self-polar triangle $DEF$. Let $$\mathcal{S}:\mathbf{P}\backslash\left\{D,E,F\right\}\rightarrow\mathbf{P}\textrm{  ,  }P\mapsto\mathcal{S}(P)=P'$$ be the Steiner correspondence with respect to $\Sigma_1$ and $\Sigma_2$. Then the sets $$\left\{U',V',W'\right\}\textrm{  ,  }\left\{U',A',C'\right\}\textrm{  ,  }\left\{V',B',C'\right\}\textrm{  ,  }\left\{W',A',B'\right\}$$ are collinear. So by observation (5), $$X'=\egy{DA'}\cap\egy{U'V'}\textrm{  ,  }Y'=\egy{DB'}\cap\egy{U'V'}\textrm{  ,  }Z'=\egy{DC'}\cap\egy{U'V'}.$$ Thus the opposite side lines of the complete quadrangle $A'B'C'D$ meet the line $\egy{U'V'}$ at the pairs of points $(U'X')$, $(V'Y')$, $(W'Z')$. These are the pairs of an involution on the line $\egy{U'V'}$, therefore, in view of (4), their images $(UX)$, $(VY)$, $(WZ)$ under $\mathcal{S}$ are the pairs of an involution on the conic $\Sigma$. By our construction this is equivalent to the property that the lines $\egy{UX}$, $\egy{VY}$, $\egy{WZ}$ are concurrent.

\footnotesize
\textsc{\\Zoltán Szilasi\\Institute of Mathematics\\University of Debrecen\\H-4010 Debrecen\\Hungary}\\
\textit{E-mail}: szilasi.zoltan@science.unideb.hu
	

\begin{thebibliography}{9}
	\bibitem{BradSmith} C. J. Bradley, G. C. Smith: On a construction of Hagge, \textit{Forum Geometricorum} \textbf{7} (2007), 231-247.
	\bibitem{B} C. J. Bradley: \textit{Generalizations of Hagge's Theorems}, arXiv:1007.2762, 2010.
	\bibitem{Cox1} H. S. M. Coxeter: \textit{The Real Projective Plane}, Second edition, Cambridge, 1955.
	\bibitem{Hagge} K. Hagge: Der Fuhrmannsche Kreis und der Brocardsche Kreis als Sonderf\"{a}lle eines allgemeineren
Kreises, \textit{Zeitschrift für Math. Unterricht} \textbf{38} (1907), 257-269.
	\bibitem{Hatt} J. L. S. Hatton: \textit{The Principles of Projective Geometry Applied to the Straight Line and Conic}, Cambridge, 1913.
	\bibitem{J} R. A. Johnson: \textit{Advanced Euclidean Geometry}, Dover, 1960.
	\bibitem{Preiser} A. M. Peiser: The Hagge circle of a triangle, \textit{American Mathematical Monthly} \textbf{49} (1942), 524-527.
\end{thebibliography}
\end{document}